\documentclass[a4paper]{gtart}

\usepackage{psfrag, graphicx, subfigure, epsfig,color}

\usepackage{amsmath, amssymb, latexsym, euscript}

\usepackage{mathptmx, wrapfig}

\usepackage[inner=30mm, outer=30mm, textheight=225mm]{geometry}

\def\lst{\mathcal{S}}
\def\tri{\mathcal{T}}

\def\even{\mathfrak{e}}

\def\bkR{{\rm I\kern-.17em R}}
\def\R{\bkR}
\def\bkZ{{\rm Z\kern-.26em Z}}
\def\Z{\bkZ}

\def\bkbZ{{\rm Z\kern-.27em Z}}
\def\bZ{\bkbZ}

\DeclareMathOperator{\e}{e}
\DeclareMathOperator{\T}{T}

\theoremstyle{plain}
\newtheorem{theorem}{Theorem}
\newtheorem*{theorem*}{Theorem}
\newtheorem{lemma}[theorem]{Lemma}
\newtheorem{proposition}[theorem]{Proposition}

\theoremstyle{definition}
\newtheorem{definition}[theorem]{Definition}
\newtheorem*{definition*}{Definition}

\theoremstyle{remark}
\newtheorem{remark}[theorem]{Remark}

\numberwithin{equation}{section}

\begin{document}

\title{$\mathbf{\bZ_2}$--Thurston Norm and Complexity of 3--Manifolds}
\author{William Jaco, J.\thinspace Hyam Rubinstein and Stephan Tillmann}

\begin{abstract}
A new lower bound on the complexity of a 3--manifold is given using the $\Z_2$--Thurston norm. This bound is shown to be sharp, and the minimal triangulations realising it are characterised using normal surfaces consisting entirely of quadrilateral discs.
\end{abstract}

\makeshorttitle


\section*{Prologue}

Given a closed, irreducible 3--manifold, its complexity is the minimum number of tetrahedra in a (pseudo--simplicial) triangulation of the manifold. This number agrees with the complexity defined by Matveev~\cite{Mat1990} unless the manifold is $S^3,$ $\R P^3$ or $L(3,1).$ It follows from the definition that the complexity is known for all closed, irreducible manifolds which appear in certain computer generated censuses. In general, the question of determining the complexity of a given closed 3-manifold is difficult and one is therefore interested in finding both upper and lower bounds. Whilst an upper bound arises from the presentation of a manifold via a spine, a Heegaard splitting or a triangulation, Matveev~\cite{Mat2003} states that the problem of finding lower bounds is quite difficult. Lower bounds using homology groups or the fundamental group are given by Matveev and Pervova~\cite{MatPer2001}, and lower bounds using hyperbolic volume are given by Matveev, Petronio and Vesnin~\cite{MPV}. These bounds are only known to be sharp for a few census examples.

In \cite{JRT2} the authors found a lower bound for the complexity using covering spaces and used it to classify all manifolds realising this lower bound. In particular, this determined two infinite families of minimal triangulations, and hence the complexity for infinitely many manifolds. 
This lower bound supposes the existence of a non-trivial $\Z_2$--cohomology class (or, equivalently, the existence of a connected double cover of the manifold). The present paper uses the existence of multiple $\Z_2$--cohomology classes to give a new lower bound for complexity using an analogue of Thurston's norm, and the minimal triangulations realising this bound are characterised. Moreover, it is shown that an infinite family of triangulations realises this bound. Whilst this family already arose in \cite{JRT2}, the methods are more generally applicable and lead to a structure theory for the minimal triangulations which are close to realising the lower bound. The new lower bound also gives very tight two sided bounds for many new examples. Moreover, the bootstrapping method of \cite{JRT2} using covers can be combined with this approach --- using higher covering degree and the existence of covers with arbitrarily large $\Z_2$--cohomology.

The effectiveness of the new bounds arising in this work result from the desire not only to know the complexity of a manifold but also some or all minimal triangulations realising this complexity. The combinatorial structure of a minimal triangulation is governed by 0--efficiency \cite{JR} and low degree edges \cite{JRT}. From this one can extrapolate building blocks for minimal triangulations. Understanding how they fit together under extra constraints on the manifold is a guiding principle in this work. Using these ideas, one can effectively try to understand vertical sections of the census which pick up a finite cover of every manifold. Especially with view towards infinite families of minimal triangulations of hyperbolic manifolds, this seems to be the most promising approach to date as the largest known census at the time of writing goes up to complexity 12, but hyperbolic examples only appear from complexity 9 and are sparse amongst these low complexity manifolds.


\section{Definitions, results and applications}

Let $M$ be a closed, orientable, irreducible, connected 3--manifold. Let $\varphi \in H^1(M;\Z_2),$ and $S$ be a properly embedded surface dual to $\varphi.$ An analogue of Thurston's norm \cite{Th86} can be defined as follows. If $S$ is connected, let $\chi_{-}(S) = \max \{ 0,-\chi(S)\},$ and otherwise let
$$\chi_{-}(S) = \sum_{S_i\subset S} \max \{ 0,-\chi(S_i)\},$$
where the sum is taken over all connected components of $S.$ Note that $S_i$ is not necessarily orientable. Define:
$$|| \ \varphi \ || = \min \{ \chi_{-}(S) \mid S \text{ dual to } \varphi\}.$$
The surface $S$ dual to $\varphi \in H^1(M;\Z_2)$ is said to be \emph{$\Z_2$--taut} if no component of $S$ is a sphere and $\chi(S) = -|| \ \varphi \ ||.$ As in \cite{Th86}, one observes that every component of a $\Z_2$--taut surface is non-separating and geometrically incompressible.

\begin{theorem}[Thurston norm bounds complexity]\label{thm:general bound}
Let $M$ be a closed, orientable, irreducible, connected 3--manifold with triangulation $\tri,$ 
and denote by $|\tri|$ the number of tetrahedra. If $H \le H^1(M; \Z_2)$ is a subgroup of rank two, then:
$$
|\tri| \ge 2 + \sum_{0 \neq \varphi \in H} || \ \varphi \ ||.
$$
\end{theorem}

There is a nice characterisation for triangulations realising the above lower bound. Let $\tri$ be a triangulation of $M$ having a single vertex. Place three quadrilateral discs in each tetrahedron, one of each type, such that the result is a (possibly branched immersed) normal surface. This surface is denoted $Q$ and called the \emph{canonical quadrilateral surface}. Suppose $Q$ is the union of three embedded normal surfaces. Then each of them meets each tetrahedron in a single quadrilateral disc and is hence a one-sided Heegaard splitting surface. It defines a dual $\Z_2$--cohomology class and $H^1(M; \Z_2)$ has rank at least two.

\begin{theorem}\label{thm:really main}
Let $M$ be a closed, orientable, irreducible, connected 3--manifold with triangulation $\tri.$ Let $H \le H^1(M; \Z_2)$ be a subgroup of rank two. Then the following two statements are equivalent.
\begin{enumerate}
\item We have
$$
|\tri| = 2 + \sum_{0 \neq \varphi \in H} || \ \varphi \ ||.
$$
\item
The triangulation has a single vertex and the canonical quadrilateral surface is the union of three $\Z_2$--taut surfaces representing the non-trivial elements of $H.$
\end{enumerate}
Note that (1) implies that $\tri$ is minimal by Theorem \ref{thm:general bound}. Moreover, (2) implies that each non-trivial element of $H$ has a $\Z_2$--taut representative, which is a one-sided Heegaard splitting surface, and that each edge has even degree. 
\end{theorem}

The proofs of the theorems are based on a refinement of the methods of \cite{JRT}. New results are given concerning combinatorial constraints from  $\Z_2$--cohomology classes and intersections of maximal layered solid tori. These results can be found in Sections~\ref{sec:co-homology classes} and \ref{sec:Intersections of maximal layered solid tori} respectively. In Section~\ref{sec:quad surfaces and 1-sided Hs}, we study quadrilateral surfaces and their relationship with Heegaard splittings. It is interesting to note that one-sided splittings of lowest complexity are obtained from the minimal triangulations in \cite{JRT} and \cite{JRT2}. In contrast, determining two-sided Heegaard splittings of lowest complexity is in general very difficult. The proofs of the main results are given in Section~\ref{sec:degree three}.

In Section~\ref{sec:triangulations}, we show that the \emph{twisted layered loop triangulation} of $S^3/Q_{8k},$ $k$ any positive integer, satisfies the equivalent statements in Theorem \ref{thm:really main}. It was already shown in \cite{JRT2} that this triangulation is the unique minimal triangulation.

We conclude this introduction with a few open problems. 

(1) Are there more manifolds (in particular with $\Z_2$--cohomology of rank three or more) satisfying the equivalent statements of Theorem~\ref{thm:really main}?

(2) Are there triangulations with more than one vertex and such that the canonical quadrilateral surface is the union of three $\Z_2$--taut surfaces?

(3) The small Seifert fibred space 
$$M_{m,n} = S^2( (1,-1), (2,1), (2m+2,1), (2n+2,1)),$$ 
where $m$ and $n$ are positive integers, is triangulated by a \emph{layered chain pair} having $2(m+n)+2$ tetrahedra (see \cite{bab}). This satisfies $ |\tri |  = 4+ \sum || \ \varphi \ ||.$ Theorem \ref{thm:really main} therefore implies that $M_{m,n}$ has complexity $2(m+n),$ $2(m+n)+1$ or $2(m+n)+2.$ What is the complexity of $M_{m,n}$?

(4) The Seifert fibred space
$$M_{k,m,n} = S^2( (1,-1), (2k+2,1), (2m+2,1), (2n+2,1)),$$ 
where $k, m$ and $n$ are positive integers, is triangulated by an \emph{augmented solid torus} having $2k+2m+2n+3$ tetrahedra (see \cite{bab}). Is this a minimal triangulation satisfying $ |\tri |  = 3 + \sum || \ \varphi \ ||$?

(5) In general, the above theorem gives a very good estimate for the complexity of manifolds having triangulations satisfying $ |\tri |  = 3+ \sum || \ \varphi \ ||$ or $ |\tri |  = 4+ \sum || \ \varphi \ ||.$ Determine a complete profile of all minimal triangulations satisfying these equalities. From work of Martelli and Petronio \cite{MP}, it appears likely that the minimal triangulations for many Seifert fibred spaces with the appropriate cohomology fall into this range.

(6) Determine an effective bound for the complexity of $M$ using a rank $k$ subgroup of $H^1(M; \Z_2)$ for $k\ge 3.$

The first author is partially supported by NSF Grant DMS-0505609 and the Grayce B. Kerr Foundation. The second and third authors are partially supported under the Australian Research Council's Discovery funding scheme (project number DP0664276).


\section{Normal surfaces dual to $\mathbf{\Z_2}$--cohomology classes}
\label{sec:co-homology classes}

Throughout this section, let $\tri$ be an arbitrary 1--vertex triangulation of the closed, orientable, connected 3--manifold $M.$ A non-trivial class in $H^1(M, \Z_2)$ was used in \cite{JRT} to study $\tri.$ This naturally generalises to subgroups of $H^1(M, \Z_2)$ of arbitrary rank. For the purpose of this paper, it suffices to consider rank--2 subgroups. To fix notation, assume that $\varphi_1, \varphi_2 \in H^1(M, \Z_2)$ such that
$$H = \langle \varphi_1, \varphi_2 \rangle \cong \Z_2 \oplus \Z_2.$$
We let $\varphi_3=\varphi_1 + \varphi_2.$ A colouring of the edges arising from $H$ is introduced and a canonical surface is associated to $H.$ It is shown that this yields a combinatorial constraint for the triangulation.


\subsection{Triangulations}

The notation of \cite{JR, JR:LT} will be used in this paper. Hence $\tri$ consists of a union of pairwise disjoint 3--simplices, $\widetilde{\Delta},$ a set of face pairings, $\Phi,$ and a natural quotient map $p\co \widetilde{\Delta} \to \widetilde{\Delta} / \Phi = M.$ Since the quotient map is injective on the interior of each 3--simplex, we will refer to the image of a 3--simplex in $M$ as a \emph{tetrahedron} and to its faces, edges and vertices with respect to the pre-image. Similarly for images of 2-- and 1--simplices, which will be referred to as \emph{faces} and \emph{edges} in $M.$ For edge $e,$ the number of pairwise distinct 1--simplices in $p^{-1}(e)$ is termed its \emph{degree}, denoted $d(e).$ If $e$ is contained in $\partial M,$ then it is a \emph{boundary edge}; otherwise it is an \emph{interior edge}.


\subsection{Rank--1 colouring of edges and canonical surface}
\label{subsec:colouring}

\begin{figure}[t]
\psfrag{1}{{\small $I$}}
\psfrag{2}{{\small $II$}}
\psfrag{3}{{\small $III$}}
\psfrag{4}{{\small $IV$}}
\psfrag{5}{{\small $V$}}
\begin{center}
    \subfigure[Rank--1 colouring]{\label{fig:Z2-homology class}
      \includegraphics[height=2.3cm]{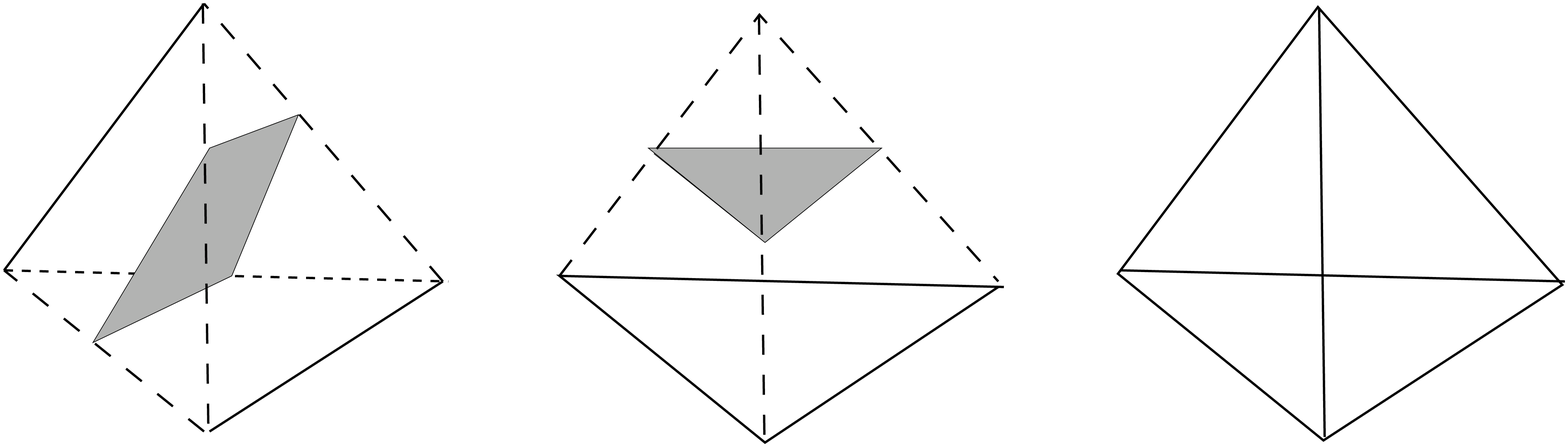}
} \\
    \subfigure[Rank--2 colouring]{\label{fig:multiZ2-homology class}
      \includegraphics[height=2.3cm]{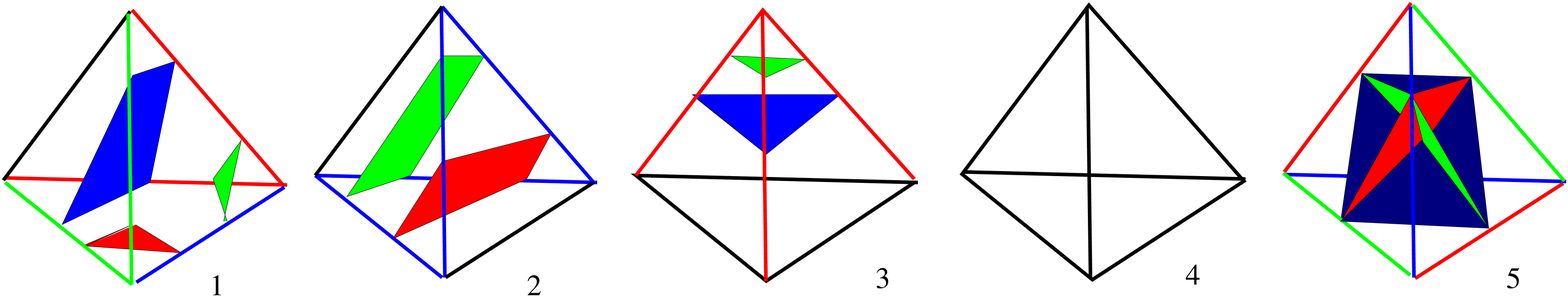}
}
\end{center}
\caption{Colouring of edges and dual normal discs}
\end{figure}

Let $0 \neq \varphi \in H^1(M, \Z_2).$ The following construction can be found in \cite{JRT}. Edge $e$ is given an orientation, and hence represents an element $[e]\in \pi_1(M).$ If $\varphi[e]=0,$ $e$ is termed $\varphi$--even, otherwise it is termed $\varphi$--odd. This terminology is independent of the chosen orientation for $e.$ Faces in the triangulation give relations between loops represented by edges. It follows that a tetrahedron falls into one of the following categories, which are illustrated in Figure~\ref{fig:Z2-homology class}:
\begin{itemize}
\item[] Type 1: A pair of opposite edges are $\varphi$--even, all others are $\varphi$--odd.
\item[] Type 2: The three edges incident to a vertex are $\varphi$--odd, all others are $\varphi$--even.
\item[] Type 3: All edges are $\varphi$--even.
\end{itemize}
If $\varphi$ is non-trivial, one obtains a unique normal surface, $S_\varphi = S_\varphi(\tri),$ with respect to $\tri$ by introducing a single vertex on each $\varphi$--odd edge. This surface is disjoint from the tetrahedra of type 3; it meets each tetrahedron of type 2 in a single triangle meeting all $\varphi$--odd edges; and each tetrahedron of type 1 in a single quadrilateral dual to the $\varphi$--even edges. Moreover, $S_\varphi$ is dual to $\varphi$ and will be termed the \emph{canonical surface dual to $\varphi.$} Since $S_\varphi$ meets each edge in the triangulation at most once, we have:

\begin{lemma}\label{lem:chi bounds norm for canonical}
Let $M$ be a closed, orientable, irreducible, connected 3--manifold and $\tri$ be a triangulation with one vertex. 
Given $0 \neq \varphi \in H^1(M, \Z_2),$ we have $$|| \ \varphi \ || \le - \chi (S_\varphi)$$ unless $M=\R P^3.$
\end{lemma}


\subsection{Rank--2 colouring of edges}
\label{subsec:rank-2 colouring}

Given the subgroup $H = \langle \varphi_1, \varphi_2 \rangle \cong \Z_2 \oplus \Z_2$ of $H^1(M; \Z_2),$ we now introduce a refinement of the above colouring. Since $\varphi_1 + \varphi_2 = \varphi_3,$ there are four types of edges:
\begin{itemize}
\item[] edge $e$ is \emph{$H$--even} or \emph{$0$--even} if $\varphi_i[e] =0$ for each $i \in \{ 1,2,3\};$ and
\item[] edge $e$ is $i$--even if $\varphi_i[e] =0$ for a unique $i \in \{ 1,2,3\}.$
\end{itemize}
Let $\{ i,j,k\} = \{ 1,2,3\}.$ The normal corners of the normal surface $S_{\varphi_i}(\tri)$ are precisely on the $j$--even and $k$--even edges. Edge $e$ is $\varphi_i$--even if it is $i$--even or $0$--even. It follows that a face of a tetrahedron either has all of its edges $0$--even; or it has two $i$--even and one $0$--even edge; or it has one $1$--even, one $2$--even and one $3$--even edge. Whence an \emph{oriented} tetrahedron falls into one of the following categories:
\begin{itemize}
\item[] Type I: One edge is $0$--even, the opposite edge is $i$--even, and one vertex of the latter is incident with two $j$--even edges, and the other with two $k$--even edges, where $\{i,j,k\} = \{ 1,2,3\}.$ (There are six distinct sub-types.)
\item[] Type II: A pair of opposite edges are $0$--even, all others are $i$--even for a unique $i \in \{ 1,2,3\}.$ (There are hence three distinct sub-types.)
\item[] Type III: The three edges incident to a vertex are $i$--even for a fixed $i \in \{ 1,2,3\},$ and all others are $0$--even. (There are hence three distinct sub-types.)
\item[] Type IV: All edges are $0$--even.
\item[] Type V: Each vertex is incident to an $i$--even edge for each $i \in \{ 1,2,3\}.$ (In particular, no edge is $0$--even, opposite edges are of the same type, and there are two distinct sub-types of tetrahedra.)
\end{itemize}
For each type, one sub-type is shown in Figure \ref{fig:multiZ2-homology class}; the black edges correspond to $0$--even edges and the normal discs in $S_{\varphi_i}$ have the same colour as the $i$--even edges. The remaining subtypes are obtained by permuting the colours other than black.


\subsection{Combinatorial bounds for triangulations}

The set-up and notation of the previous subsection is continued. Let
\begin{itemize}
\item[] $A(\tri)=$ number of tetrahedra of type I,
\item[] $B(\tri)=$ number of tetrahedra of type II,
\item[] $C(\tri)=$ number of tetrahedra of type III,
\item[] $D(\tri)=$ number of tetrahedra of type IV,
\item[] $E(\tri)=$ number of tetrahedra of type V,
\item[] $\even(\tri)=$ number of $0$--even edges,
\item[] $\tilde{\even}(\tri)=$ number of pre-images of $0$--even edges in $\widetilde{\Delta}.$
\end{itemize}

The number of tetrahedra in $\tri$ is $T(\tri)= A(\tri)+B(\tri)+C(\tri)+D(\tri)+E(\tri).$ For the remainder of this subsection, we will write $A = A(\tri),$ etc. 

\begin{lemma}
$C$ and $E$ are even.
\end{lemma}

\begin{proof}
The colouring of edges is pulled back to $\widetilde{\Delta}.$ Then the number of faces in $\widetilde{\Delta}$ having all edges $0$--even is $C+4D.$ Since faces match up in pairs and no face is identified with itself, it follows that $C$ is even.

Now consider the normal surface $S_{\varphi_1}(\tri).$ It meets an edge in a normal corner if and only if the edge is either $2$--even or $3$--even. We examine the sum of the $Q$--matching equations of all $2$--even edges (see \cite{to}). This sum, $\Sigma,$ must equal zero. Since $M$ is orientable, opposite corners of each quadrilateral disc have the same sign, and adjacent corners have opposite signs (see \cite{ti}). Let $\square$ be a quadrilateral disc in $S_{\varphi_1}(\tri).$ If $\square$ is contained in a tetrahedron of type I and meets a $2$--even edge, then it has precisely two adjacent corners on $2$--even edges; hence the contribution to $\Sigma$ of the normal corners of $\square$ is $(+1) + (-1)=0.$ If $\square$ is contained in a tetrahedron of type II, and meets a $2$--even edge, then it meets a $2$--even edge with each of its corners; hence the contribution to $\Sigma$ is $2(+1) + 2(-1)=0.$ If $\square$ is contained in a tetrahedron of type V and meets a $2$--even edge, then it meets $2$--even edges in a pair of diagonally opposite corners; hence the contribution to $\Sigma$ is either $2(+1)$ or  $2(-1).$ It follows that the number of tetrahedra of type IV, $E,$ must be even.
\end{proof}

\begin{lemma}\label{lem:char of small cases}
Let $M$ be a closed, orientable, irreducible 3--manifold with minimal triangulation $\tri.$ Suppose that all edge loops are coloured by the rank--2 subgroup $H$ of $H^1(M; \Z_2).$ Suppose that $A+C \le 3.$ 
Then $(A,B,C,D,E)$ is of one of the following forms:
$$(0,0,0,0,E), (2,B,0,0,E), (3,B,0,0,E),$$
and if $A \neq 0,$ then there is a unique edge incident with all tetrahedra of type I.
\end{lemma}

\begin{proof}
We distinguish the cases $A=0$ and $A \neq 0.$

(Case 1a) $A=0$ and $E \neq 0.$ Since each face of a tetrahedron of type II, III or IV contains an $H$--even edge and $M$ is connected, we have $B=C=D=0.$

(Case 1b) $A=0$ and $E=0.$ In this case, the Haken sum $S_{\varphi_1}+S_{\varphi_2}+S_{\varphi_3}$ is defined and is isotopic to the boundary of a regular neighbourhood, $N,$ of the complex $K$ spanned by all $0$--even edges. Then $M\setminus N$ therefore meets a tetrahedron either in the empty set or in a product region and hence each component of $M\setminus N$ is an $I$--bundle or a twisted $I$--bundle over a surface. If $C=0,$ then either all tetrahedra are of type II or all tetrahedra are of type IV. In each case, this contradicts the fact that $H^1(M; \Z_2)$ has rank two. Hence $C=2$ since $C$ is even. Moreover, the two tetrahedra of type III must be of distinct sub-types since otherwise $H^1(M; \Z_2)$ has rank one. We may assume that one of them has three $0$--even and three $1$--even edges. Pulling back the colouring of faces to $\widetilde{\Delta},$ we have an odd number of faces with one $0$--even edge and two $1$--even edges. This contradicts the fact that $M$ is closed. Hence, $A=0$ and $E=0$    
can not both occur.

(Case 2) $A\neq 0.$ The abstract neighbourhood of each 0--even edge contained in a tetrahedron of type I contains either at least two distinct tetrahedra of type I, or one tetrahedron of type one and at least two tetrahedra of type III. Since $C$ is even and $A+C \le 3,$ we have the following cases: 
($A\in \{1,2, 3\}, C=0$) and $(A=1, C=2).$ 

We first show that $(A=1, C=2)$ is not possible. Note that there is a unique 0--even edge, $e_0,$ incident with the tetrahedron of type I, $\sigma_0.$ It follows that it must be incident with the two tetrahedra of type III, $\sigma_1$ and $\sigma_2.$ Denote the two remaining $0$--even edges incident with $\sigma_1$ by $e_1$ and $e_2,$ oriented such that $\sigma_0 + \sigma_1 + \sigma_2$ is homologically the boundary of the face. Consider the abstract neighbourhood $B(e_1)$ of $e_1,$ and recall that no face can be a cone or a dunce hat. It follows that either $\sigma_0 \neq \sigma_1=\sigma_2$ or $\sigma_0 = \sigma_1=\sigma_2.$ Give $B(e_1)$ an orientation. Then  two of the tetrahedra mapping to $\sigma_1$ under the map $B(e_1) \to M$ induce opposite orientations on $\sigma_1,$ contradicting the fact that $M$ is orientable.

It is easy to see that $(A=1, C=0)$ is not possible by examining the neighbourhood of the unique 0--even edge, $e_0,$ incident with the tetrahedron of type I.

The remaining cases are $(A=2, C=0)$ and $(A=3, C=0).$ In each case, there is a unique 0--even edge, $e_0,$ contained in all tetrahedra of type I. Note that $e_0$ cannot be of degree three since otherwise the classification of degree three edges in \cite{JRT} implies that $e_0$ is contained in a maximal layered torus subcomplex and hence every tetrahedron incident with it it is of type II or IV. Moreover, if $A \neq 0$ but $C=0,$ then necessarily $D=0$ since $M$ is connected.
\end{proof}

Let $K$ be the complex in $M$ spanned by all $0$--even edges. Let $N$ be a small regular neighbourhood of $K.$ Then $\partial N$ is a normal surface; it meets each tetrahedron in the same number and types of normal discs as $S_{\varphi_1} \cup S_{\varphi_2}\cup S_{\varphi_3}$ except for the tetrahedra of type V, which it meets in four distinct normal triangle types instead of three distinct normal quadrilateral types. The normal coordinate of $\partial N$ is thus obtained by taking the sum of the normal coordinates of $S_{\varphi_1},$ $S_{\varphi_2},$ and $S_{\varphi_3},$ and adding to this the \emph{tetrahedral solution} of each tetrahedron of type V. (The tetrahedral solution is obtained by adding all triangle coordinates of a tetrahedron and subtracting all quadrilateral coordinates; see, for instance, \cite{kr}.) Hence
$$
\chi(S_{\varphi_1}) +\chi (S_{\varphi_2}) + \chi (S_{\varphi_3}) + E = \chi (\partial N) = 2 \chi(N) = 2 \chi(K) = 2 - 2 \even +C + 2D,$$
where the Euler characteristic of $K$ is computed from the combinatorial data. In particular, 
$\chi(S_{\varphi_1}) +\chi (S_{\varphi_2}) + \chi (S_{\varphi_3}) \text{  is even.}$
Rearranging the above equality gives
\begin{equation}\label{eq:e}
C + 2D -E = 2 \even - 2 + \chi(S_{\varphi_1}) +\chi (S_{\varphi_2}) + \chi (S_{\varphi_3}),
\end{equation}
and hence:
\begin{align}\label{inequ for one vert tri gen}
\tilde{\even} 	&= A+2B+3C+6D\nonumber \\
		&=  2T  -A + C +4D-2E  \nonumber\\
		&= 2T  - A - C+4 \even -4+ 2 (\chi(S_{\varphi_1}) +\chi (S_{\varphi_2}) + \chi (S_{\varphi_3})).
\end{align}
\begin{lemma}\label{lem:formula for degree 3 edges}
Let $M$ be a closed, orientable, irreducible 3--manifold, and suppose that $\varphi_1, \varphi_2 \in  H^1(M; \Z_2)$ are non--trivial with $\varphi_1+\varphi_2 = \varphi_3 \neq 0.$ Let $\tri$ be a minimal triangulation with $T$ tetrahedra, and let $S_{\varphi_i}$ be the canonical surface dual to $\varphi_i.$ Letting $\even_d$ denote the number of $0$--even edges of degree $d,$ we have:
\begin{equation}\label{inequ for min tri general}
\even_3= 4 +A +C - 2 (T + \chi(S_{\varphi_1}) +\chi (S_{\varphi_2}) +\chi (S_{\varphi_3})) + \sum_{d=5}^{\infty} (d-4) \even_d.
\end{equation}
\end{lemma}

\begin{proof}
First note that the existence of $\varphi_1$ and $\varphi_2$ implies that $M$ is not homeomorphic to one of $S^3,$ $L(3,1),$ $\R P^3,$ or $L(4,1).$ It now follows from \cite{JR}, Theorem 6.1, that $\tri$ has a single vertex; hence $S_{\varphi_i}$ and $\even_i$ are defined. Moreover, \cite{JR}, Proposition 6.3 (see also Proposition 8 of \cite{JRT}), implies that the smallest degree of an edge in $\tri$ is three. One has $\tilde{\even} = \sum d \even_d$ and $\even = \sum \even_d.$ Putting this into (\ref{inequ for one vert tri gen}) gives the desired equation.
\end{proof}


\section{Intersections of maximal layered solid tori}
\label{sec:Intersections of maximal layered solid tori}

Throughout this section, let $M$ be a closed, irreducible, orientable, connected 3--manifold with triangulation $\tri.$ We extend results of \cite{JRT} concerning maximal layered solid tori in $M.$


\begin{figure}[t]
\psfrag{M}{{\small $N$}}
\psfrag{S}{{\small $\Delta$}}
\psfrag{e}{{\small $e_1$}}
\psfrag{f}{{\small $e_2$}}
\psfrag{g}{{\small $e_3$}}
\begin{center}
        \subfigure[The solid torus $\lst_1$]{\label{fig:solid_torus}
      \includegraphics[height=3.6cm]{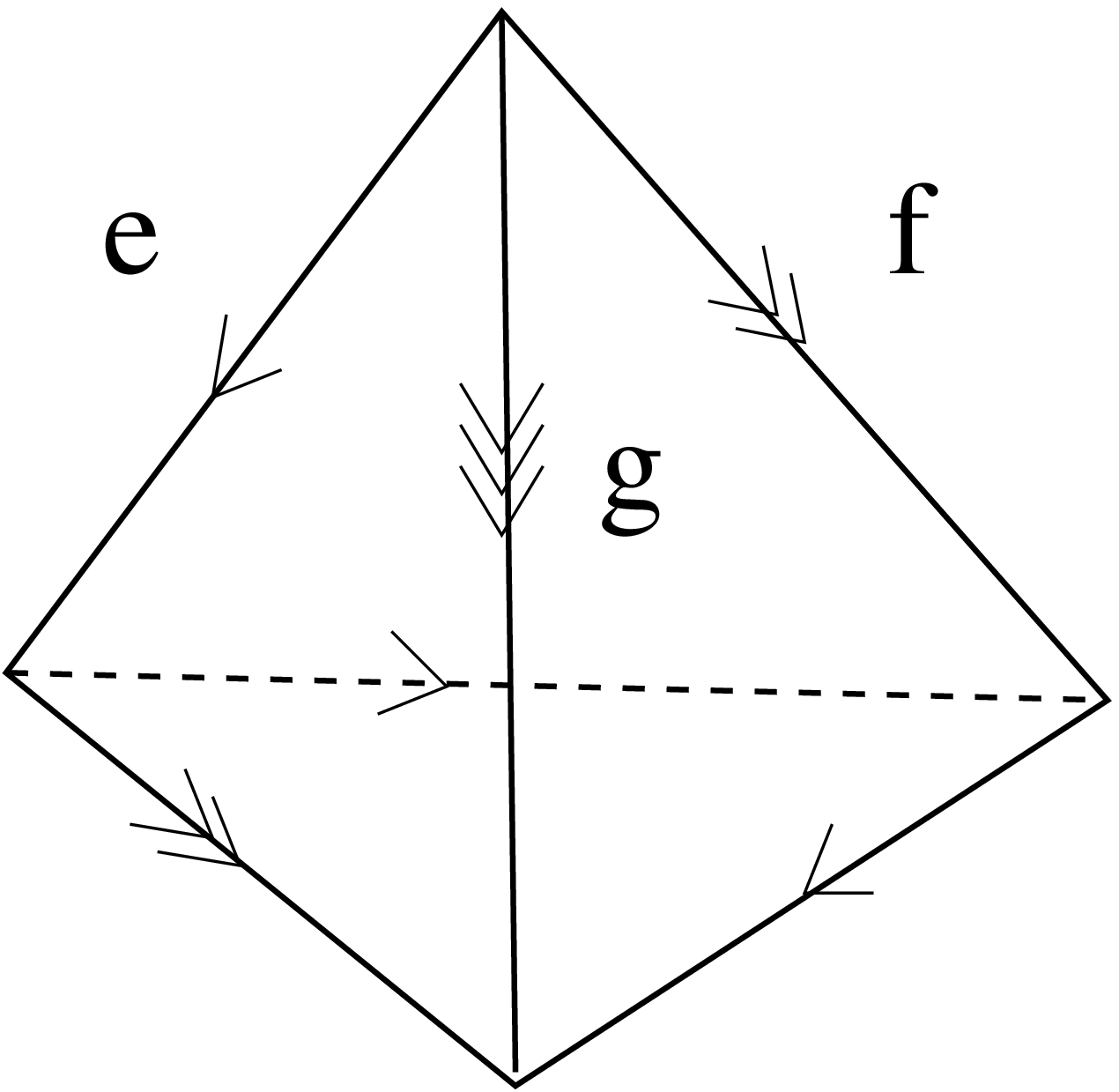}
    }
    \qquad\qquad 
        \subfigure[Layering along a boundary edge]{\label{fig:layering}
      \includegraphics[height=3.6cm]{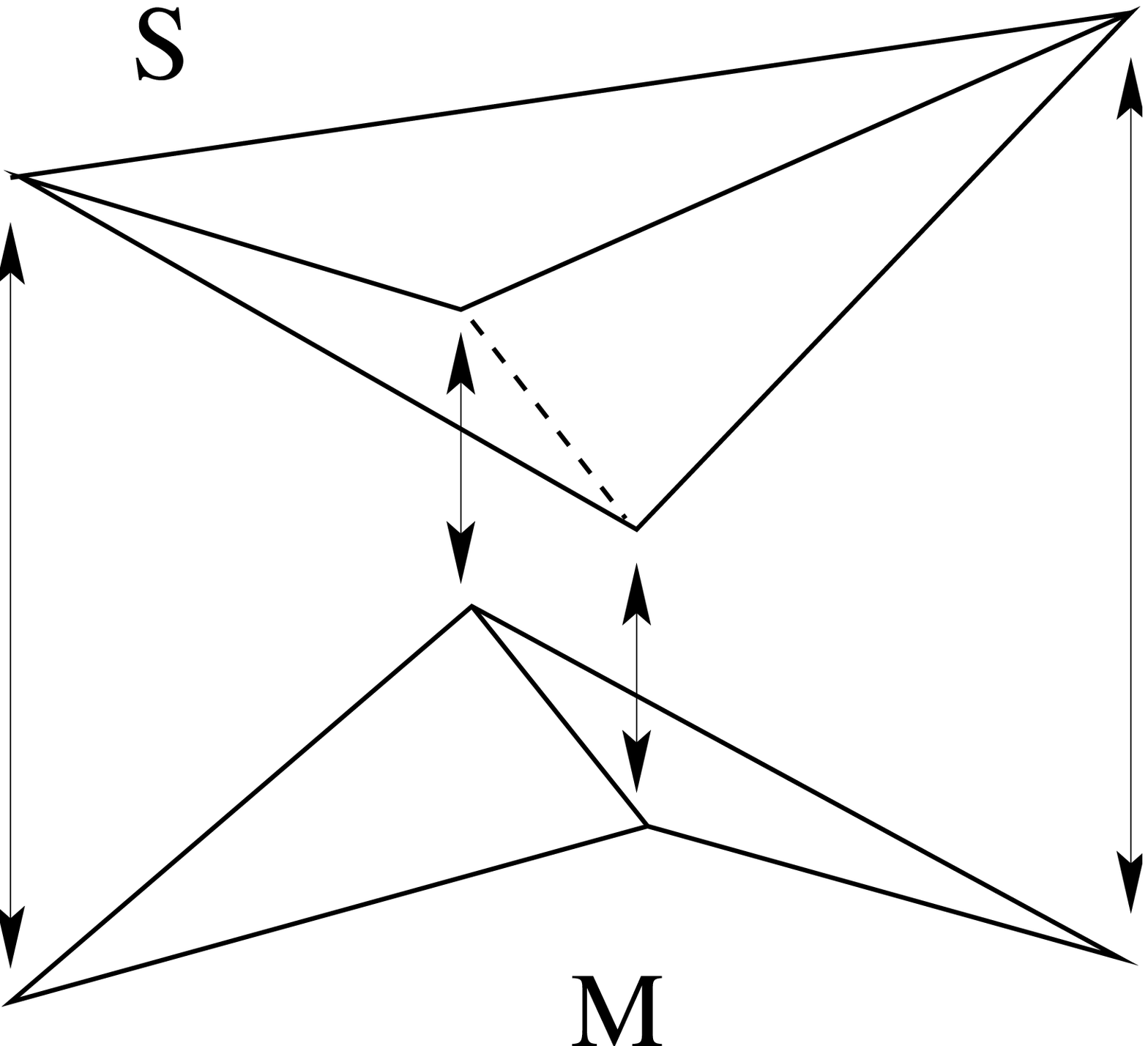}
    } 
\end{center}
    \caption{Layered triangulation of the solid torus}
     \label{fig:regular exchange is associative}
\end{figure}

\subsection{Layered solid tori}

The following definitions and facts can be found with more detail in \cite{JR:LT} and \cite{JRT}. A \emph{layered solid torus} is a solid torus with a special triangulation: the triangulation is obtained from the triangulation given in Figure \ref{fig:solid_torus} by iteratively \emph{layering along boundary edges}. Layering along a boundary edge is illustrated in Figure \ref{fig:layering}. Namely, suppose $N$ is a 3--manifold, $\tri_\partial$ is a triangulation of $\partial N$, and $e$ is an edge in $\tri_\partial$ which is incident to two distinct faces. We say the 3--simplex $\sigma$ is \emph{layered along} $e$ if two faces of $\sigma$ are paired, ``without a twist," with the two faces of $\tri_\partial$ incident with $e$. The resulting 3--manifold is homeomorphic with $N$. If $\tri_\partial$ is the restriction of a triangulation of $N$ to $\partial N$, then we get a new triangulation of $N$ and a new triangulation of $\partial N$ which differs from $\tri_\partial$ by a diagonal flip. For every edge in $M,$ we will also refer to its degree as its $M$--degree, and its degree with respect to the layered solid torus $T$ in $M$ is called its $T$--degree.

\begin{definition}(Layered solid torus in $M$)
A \emph{layered solid torus with respect to $\tri$} in $M$ is a subcomplex in $M$ which is combinatorially equivalent to a layered solid torus. 
\end{definition}

\begin{definition}(Maximal layered solid torus in $M$)
A layered solid torus is a \emph{maximal layered solid torus with respect to $\tri$} in $M$ if it is not strictly contained in any other layered solid torus in $M.$
\end{definition}

\begin{lemma}\label{lem:intersection of maximal layered solid tori}\cite{JRT}
Assume that the triangulation is minimal and 0--efficient. If $M$ is not a lens space with layered triangulation, then the intersection of two distinct maximal layered solid tori in $M$ consists of at most a single edge.
\end{lemma}

\begin{lemma}\label{lem:layered lens char 2}\cite{JRT}
Assume that the triangulation is minimal and 0--efficient, and suppose that $M$ contains a layered solid torus, $T,$ made up of at least two tetrahedra and having a boundary edge, $e,$ which has degree four in $M.$ Then either
\begin{enumerate}
\item $T$ is not a maximal layered solid torus in $M;$ or
\item $e$ is the univalent edge for $T$ and it is contained in four distinct tetrahedra in $M;$ or 
\item $M$ is a lens space with minimal layered triangulation.
\end{enumerate}
\end{lemma}

Recall the notion of rank--2 colouring from Subsection \ref{subsec:rank-2 colouring}.

\begin{lemma}
Assume that the triangulation contains a single vertex and that all edge loops are coloured by the rank--2 subgroup $H$ of $H^1(M; \Z_2).$ Then all tetrahedra in a layered solid torus in $M$ are either of type II or type IV, but not both.
\end{lemma}
\begin{proof}
The colouring of the layered solid torus is uniquely determined by the image of the longitude under the elements of $H;$ the result now follows from the layering procedure.
\end{proof}

\begin{definition}[(Types of layered solid tori)]
Assume that the triangulation contains a single vertex and that all edge loops are coloured by the rank--2 subgroup $H$ of $H^1(M; \Z_2).$ A layered solid torus containing a tetrahedron of type II (respectively IV) is accordingly termed of type II (respectively IV). 
\end{definition}


\subsection{Maximal layered solid tori in atoroidal manifolds}

An orientable 3--manifold is termed \emph{atoroidal} if it does not contain an embedded, incompressible torus.

\begin{lemma}\label{lem:normal torus in 0-eff bounds solid torus}\cite{JRT}
Assume that $\tri$ is minimal and 0--efficient and that $M$ is atoroidal. Then every torus which is normal with respect to $\tri$ bounds a solid torus in $M$ on at least one side.
\end{lemma}

Note that if there is a rank--2 subgroup $H$ of $H^1(M; \Z_2),$ then every minimal triangulation of $M$ is 0--efficient.

\begin{lemma}\label{lem:three meet in edge}
Assume that $\tri$ is minimal and that $M$ is atoroidal. Suppose the edges are coloured by the rank--2 subgroup $H$ of $H^1(M; \Z_2),$ and that three pairwise distinct maximal layered solid tori of type II meet in a $H$--even edge. Then these are the only maximal layered solid tori of type II in the triangulation and $M$ is a Seifert fibred space with base $S^2$ and precisely three exceptional fibres.
\end{lemma}

\begin{proof}
Denote the three maximal layered solid tori by $T_1, T_2, T_3,$ and the common $H$--even edge by $e.$ If $T_i \cap T_j$ properly contains $e$ for $i\neq j,$ then $M$ is a lens space with layered triangulation according to Lemma \ref{lem:intersection of maximal layered solid tori}. But this contradicts the assumption that $H^1(M; \Z_2)$ has rank at least two. Hence $T_i \cap T_j=\{ e\}.$

Let $N$ be a small regular neighbourhood of $T_1\cup T_2 \cup T_3.$ Then $\partial N$ is a (topological) torus and a barrier surface. Hence, either $\partial N$ is isotopic to a normal surface or $\overline{M\setminus N}$ is a solid torus. In the first case, Lemma \ref{lem:normal torus in 0-eff bounds solid torus} implies that either $N$ or $\overline{M\setminus N}$ is a solid torus. However, $N$ cannot be a solid torus as $e$ is not a longitude of either $T_1,$ $T_2$ or $T_3.$ Hence $\overline{M\setminus N}$ is a solid torus and, in particular, $N$ admits a Seifert fibration with three exceptional fibres, and $M$ admits a Seifert fibration with three or four exceptional fibres. If $M$ admits a Seifert fibration with four exceptional fibres, then it contains an embedded, incompressible (vertical) torus, contradicting the assumption that $M$ is atoriodal. Similarly, if $M$ admits a Seifert fibration with three exceptional fibres but the base is not $S^2,$ then $M$ contains an embedded, incompressible (vertical) torus. Hence $M$ is a Seifert fibred space with base $S^2$ and precisely three exceptional fibres.

It follows that $e$ is homotopic to a longitude of $\overline{M\setminus N}$ and each $\varphi \in H$ restricted to $\overline{M\setminus N}$ is trivial. Suppose $T_4$ is a maximal layered solid torus of type II which is distinct from $T_1, T_2, T_3.$ Then the longitude of $T_4$ is not $H$--even. But it is clearly homotopic into $\overline{M\setminus N},$ contradicting the fact that each $\varphi$ restricted to $\overline{M\setminus N}$ is trivial. Whence $T_1,$ $T_2$ and $T_3$ are the only maximal layered solid tori of type II in the triangulation.
\end{proof}

\begin{lemma}\label{lem:two meet in edge}
Assume that $\tri$ is minimal and that $M$ is atoroidal. Suppose the edges are coloured by the rank--2 subgroup $H$ of $H^1(M; \Z_2),$ and that precisely two distinct maximal layered solid tori of type II meet in an $H$--even edge, $e.$ Then $e$ is the only $H$--even edge in the triangulation which is contained in more than one maximal layered solid torus of type II and $M$ is a Seifert fibred space with base $S^2$ and precisely three exceptional fibres.
\end{lemma}

\begin{proof}
Denote the two maximal layered solid tori of type II meeting in $e$ by $T_1, T_2.$ Then, as above, the boundary of a small regular neighborhood $N$ of $T_1\cup T_2$ is a (topological) torus and a barrier surface, and it follows that $\overline{M\setminus N}$ is a solid torus. If $e$ is homotopic to a longitude of $\overline{M\setminus N},$ then $M$ admits a Seifert fibration with two exceptional fibres. Hence either $M$ is toroidal or $M$ is a lens space; the first is not possible and the second contradicts the assumption that $H^1(M; \Z_2)$ has rank at least two. Hence $e$ is not homotopic to a longitude of $\overline{M\setminus N}.$ Hence, $M$ is a Seifert fibred space with precisely three exceptional fibres, and the base must be $S^2$ since there are no vertical incompressible tori.

Now assume that $T'_1$ and $T'_2$ are two distinct maximal layered solid tori of type II meeting in the $H$--even edge $e' \neq e.$ Then $e'$ cannot be the longitude of either $T'_1$ nor $T'_2.$ But this contradicts the fact that the interiors of $T'_1$ and $T'_2$ are contained in the solid torus $\overline{M\setminus N}.$
\end{proof}


\section{Quadrilateral surfaces and Heegaard splittings}
\label{sec:quad surfaces and 1-sided Hs}

Let $M$ denote a closed, orientable, irreducible, connected 3--manifold throughout this section.


\subsection{Quadrilateral surfaces}

Let $\tri$ be a triangulation of $M.$ Place three quadrilateral discs in each tetrahedron, one of each type, such that the result is a (possibly branched immersed) normal surface, denoted $Q$ and called the \emph{canonical quadrilateral surface}. The surface $Q$ can be viewed as the image in $M$ of a surface $Q'$ with a cell decomposition into quadrilaterals, and the image of a connected component of $Q'$ will be termed a component of $Q.$

If $S$ is a component of $Q$ and meets some tetrahedron, $\sigma,$ in $i$ normal quadrilaterals, $i \in \{ 0,1,2,3\},$ then the same is true for each tetrahedron meeting $\sigma$ in a face. Since $M$ is connected, the same is true for every tetrahedron in the triangulation (and in particular $i \neq 0$ since $S\neq \emptyset$). It follows that $Q$ has at most three components.

Suppose that the component $S$ of $Q$ is embedded. Then $S$ meets every tetrahedron in a single quadrilateral disc. Such an embedded normal surface consisting entirely of quadrilateral discs, one in each tetrahedron, is called a \emph{quadrilateral surface}. 

Let $S$ be a quadrilateral surface, and let $C(S)$ denote the 1--complex in $M$ consisting of all edges disjoint from $S.$ Notice that $S$ divides each tetrahedron into two prisms. Then $M\setminus C(S)$ is foliated by $S$ and infinitely many copies of the boundary of a regular neighbourhood of $S.$ It follows that $M \setminus S$ is either an open handlebody or the disjoint union of two open handlebodies having the same genus, depending on whether $S$ is non-separating or separating. In particular, $S$ is either a one-sided or a two-sided Heegaard splitting surface for $M.$ A set of discs for the Heegaard splitting is dual to the set of edges in $C(S).$ The disc associated to edge $e$ in $C(S)$ is naturally triangulated with one triangle for each prism in $M \setminus S$ containing $e.$ This system of discs is possibly larger than a standard set of meridian discs if the triangulation has many vertices.

Let $e$ be an edge which the quadrilateral surface $S$ meets. Since $S$ meets each tetrahedron in a single quadrilateral disc, and $e$ is incident with as many quadrilateral discs of positive slope as of negative slope by the $Q$--matching equations \cite{to}, it follows that $e$ must have even degree. In particular, if the canonical quadrilateral surface $Q$ has three components, then each is a quadrilateral surface and it follows that every edge in the triangulation has even degree. 

It follows from this discussion that the existence of a canonical quadrilateral surface $Q$ with three components places strong constraints on the Heegaard diagram associated to each of its components. For instance, (1) the intersection of component $S$ with any of the two other components of $Q$ is a spine for $S;$ and (2) since edges dual to the other surfaces must be of even order, it follows that each meridian on $S$ contains an even number of intersections with itself and other meridians.


\subsection{Heegaard splittings}

The above construction of a Heegaard splitting from a quadrilateral surface can be reversed. Given the one-sided (resp.\thinspace two-sided) Heegaard splitting surface $S$ for $M,$ choose a complete system of meridian discs. Then $S$ together with the discs is a simple spine for the once (resp.\thinspace twice) punctured manifold. The usual dual cell decomposition gives a triangulation with one (resp.\thinspace two) vertices and one tetrahedron for each intersection point of the meridians. Moreover, $S$ sits in this triangulation as a quadrilateral surface.

It follows that the complexity of $M$ gives a lower bound for the complexity of a Heegaard splitting. In case of a two-sided Heegaard splitting, these complexities will not coincide unless $M$ has a two-vertex minimal triangulation. But then $M=S^3, \R P^3$ or $L(3,1),$ whence the complexities only coincide for $M=S^3.$ In contrast, the minimal triangulations in \cite{JRT} and \cite{JRT2} are all dual to one-sided Heegaard splittings and hence determine least complexity one-sided Heegaard splittings.


\subsection{Consequences for one-vertex triangulations}

\begin{lemma}\label{lem:quad surfaces}
Let $M$ be a closed, orientable, irreducible, connected 3--manifold with one-vertex triangulation $\tri.$
If the canonical quadrilateral surface $Q$ has three components, $Q_1, Q_2, Q_3,$ then each $Q_i$ is dual to some element $0 \neq \varphi_i \in H^1(M; \Z_2),$ and $\langle \varphi_1, \varphi_2, \varphi_3\rangle \cong \Z_2 \oplus \Z_2.$ Moreover,
$$ |\tri| + \sum \chi(Q_i) =2.$$
\end{lemma}

\begin{proof}
Let $v$ denote the single vertex. Then each component of $Q$ is non-separating. Since $\pi_1(M; v)$ is generated by the edges, an element $\varphi_i \in H^1(M; \Z_2)$ is defined by letting $\varphi_i[e] =0$ if $e$ is contained in $M\setminus Q_i,$ and $\varphi_i[e] =1$ otherwise. Hence $\varphi_i \neq 0$ and $\varphi_i + \varphi_j = \varphi_k,$ where all three indices are distinct. In particular, $H^1(M; \Z_2)$ has rank at least two. The first equation follows from a simple Euler characteristic argument.
\end{proof}


\section{Proofs of the main results}
\label{sec:degree three}

Let $M$ be a closed, orientable, irreducible, connected 3--manifold with minimal triangulation $\tri.$ Suppose that there is a rank--2 subgroup $H$ of $H^1(M; \Z_2).$ Then every minimal triangulation of $M$ is 0--efficient. In particular, there are no normal projective planes, and every normal 2--sphere is vertex linking.


\subsection{Promoting triangulations using edge flips}

The material of this subsection refines the argument in the proof of Theorem 5 in \cite{JRT}. A maximal layered solid torus is \emph{supportive} if it is of type II and each $H$--even edge incident with it has degree at most four. A maximal layered solid torus is \emph{almost supportive} if it is of type II and each $H$--even edge incident with it has $M$--degree at most four except for one $H$--even edge which has $M$--degree five. A maximal layered solid torus is of \emph{type (II,4)} if it is of type II and its $H$--even boundary edge has $M$--degree four. Every supportive maximal layered solid torus is of type (II,4). A minimal triangulation is termed \emph{(II,4)--free} if it does not contain any maximal layered solid torus of type (II,4).

An \emph{edge flip} is the replacement of an edge of degree four which is incident with four distinct tetrahedra by another edge incident with four distinct tetrahedra, see Figure \ref{fig:edge flip}. An edge flip may change the triangulation and the degree sequence associated to the edges, but it does not alter the number of tetrahedra.

\begin{lemma}\label{lem:no supportive}
Let $M$ be a closed, orientable, irreducible 3--manifold with minimal triangulation $\tri.$ Suppose that all edge loops are coloured by the rank--2 subgroup $H$ of $H^1(M; \Z_2).$ Then there is a minimal triangulation which is (II,4)--free and which is obtained from $\tri$ by a finite number of edge flips.
\end{lemma}

\begin{proof}
Let $T$ be a maximal layered solid torus of type (II,4). Lemma~\ref{lem:layered lens char 2} implies that $e$ is the univalent edge of $T,$ and that it is contained in four distinct tetrahedra in the triangulation. The argument proceeds by replacing the four tetrahedra around $e$ by a different constellation of four tetrahedra using an appropriate edge flip. The resulting triangulation is also minimal, and it is shown to either contain fewer maximal layered solid tori of type (II,4) or fewer tetrahedra of type IV. Since both of these numbers are finite, we arrive at a (II,4)--free minimal triangulation after a finite number of edge flips.

\begin{figure}[t]
\psfrag{a}{{\small $-1$}}
\psfrag{b}{{\small $0$}}
\psfrag{c}{{\small $+1$}}
\begin{center}
      \includegraphics[height=3.6cm]{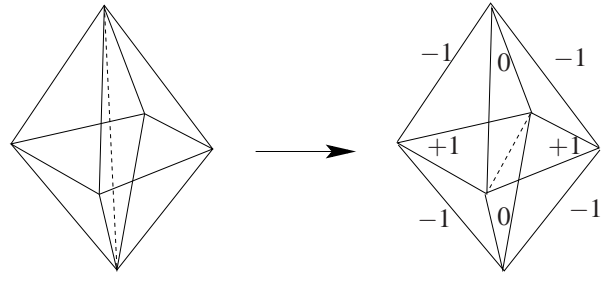}
\end{center}
    \caption{Edge flip: Labels indicate the changes in edge degree}
     \label{fig:edge flip}
\end{figure}

It remains to describe the re-triangulation process. There are various different cases to consider; they are listed below using the types of tetrahedra ordered cyclically around the edge $e,$ starting with the supportive maximal layered solid torus $T.$ The arguments given in \cite{JRT} apply almost verbatim, and we therefore merely list the case in \cite{JRT}, the corresponding neighbourhoods and the appropriate flip.

(1,1,1,1) in \cite{JRT}: Corresponds to (II,II,II,II), (II,II,I,I), (II,I,II,I), (II, I,I,II), (II,I,I,I); any edge flip decreases the total number of maximal layered solid tori of type (II,4).

(1,2,2,1) in \cite{JRT}: Corresponds to (II,III,III,II), (II, III, III, I); the edge flip with the property that the $H$--even edges (other than $e$) contained in the faces of the type II tetrahedra remain of the same degree decreases the total number of maximal layered solid tori of type (II,4).

(1,1,2,2) in \cite{JRT}: Corresponds to (II,II,III,III), (II, I, III, III); the edge flip with the property that the $H$--even edges (other than $e$) contained in the faces of the type II tetrahedra remain of the same degree decreases the total number of maximal layered solid tori of type (II,4).

(1,2,3,2) in \cite{JRT}: Corresponds to (II, III, IV, III); any edge flip reduces the number of tetrahedra of type IV.
\end{proof}


\subsection{Controlling triangulations using degree three edges}

The proof of the following lemma contains a more explicit analysis which is not contained in the statement.

\begin{lemma}\label{lem:bound for degree 3 edges}
Let $M$ be a closed, orientable, irreducible 3--manifold with minimal triangulation $\tri.$ Suppose that all edge loops are coloured by the rank--2 subgroup $H$ of $H^1(M; \Z_2)$ and that $\tri$ is (II,4)--free.
Letting $\even_d$ denote the number of $H$--even edges of degree $d,$ we have:
\begin{equation}\label{eqn:inequality degree three}
\even_3 < 2+ \sum_{d=5}^{\infty} (d-4) \even_d.
\end{equation}
\end{lemma}

\begin{proof}
Assume, by way of contradiction, that
\begin{equation}\label{eqn:degree three proof}
\even_3 \ge 2 + \sum_{d=5}^{\infty} (d-4) \even_d.
\end{equation}
The argument refines the counting argument in the proof of Theorem 5 in \cite{JRT}. Each edge of degree three is $H$--even, so (\ref{eqn:degree three proof}) implies that there are at least two edges of degree three. Since the triangulation contains at least three tetrahedra, it follows from Proposition 9 in \cite{JRT} that each edge of degree three, $e,$ is the base edge of a layered solid torus subcomplex isomorphic to $\lst_2=\{4,3,1\}.$ This subcomplex is contained in a unique maximal layered solid torus, $\T(e).$ Conversely, if a maximal layered solid torus, $T,$ contains an edge, $e,$ of degree three, then $e$ is unique and we write $e=\e(T);$ otherwise, we let $\e(T) = \emptyset.$

It follows from Lemma~\ref{lem:intersection of maximal layered solid tori} that any two distinct maximal layered solid tori share at most an edge. We seek a contradiction guided by inequality (\ref{eqn:degree three proof}). The proof, a basic counting argument, is organised as follows. The set of all edges of degree three, $Y,$ is divided into pairwise disjoint subsets; denote the set of these subsets $S(Y).$ The set of all $H$--even edges, $X,$ is similarly divided into pairwise disjoint subsets, giving a set $S(X).$ In defining these subsets, we also define an injective map $S(Y)\to S(X).$ If $Y_i \in S(Y)$ is associated with $X_i \in S(X),$ then the quantity 
$$|Y_i| - \sum_{e\in X_i} (d(e)-4)$$
is termed a \emph{deficit} if it is negative, a \emph{cancellation} if it is zero, and a \emph{gain} if it is positive. Then (\ref{eqn:degree three proof}) implies that the total gain is at least two.

Let $e \in Y$ such that $\T(e)$ is of type IV. Then $\lst_2 \cong T_0 \subseteq \T(e).$ Denote $e_0$ the longitude of $T_0.$ This is a $0$--even boundary edge with $T_0$--degree $5$ and $M$--degree $5+m$ for some $m\ge1.$ The total number of maximal layered solid tori in $M$ meeting in $e_0$ is bounded above by $\frac{m+1}{2},$ since no two meet in more than an edge. Hence the maximal layered solid tori containing a degree three edge and meeting in the $0$--even edge $e_0$ contribute at most $\frac{m+1}{2}$ to the left hand side of (\ref{eqn:degree three proof}). The contribution of $e_0$ to the right hand side is $d(e_0)-4=1+m.$ One therefore obtains a \emph{deficit} less or equal to $-\frac{m+1}{2}$ To $X_0=\{ e_0\}$ associate the set, $Y_0,$ of all degree three edges $e'$ such that $\T(e')$ contains $e_0.$

We now proceed inductively. Let $e$ be an edge of degree three such that $\T(e)$ is of type IV and $e$ is not contained in the collection of subsets $Y_0, ..., Y_{i-1}$ of $Y.$ Then $\T(e)$ contains a subcomplex isomorphic to $\lst_{2},$ whose longitude, $e_i,$ cannot be any of the edges $e_0,...,e_{i-1}.$ Consider the set $Y_i$ of all degree three edges $e'$ such that $\T(e')$ contains $e_i$ and $e'$ is not contained in any of $Y_0, ..., Y_{i-1}.$ Then the above calculation shows that there is a deficit associated to $X_i=\{e_i\}$ and $Y_i.$

It follows that there must also be a maximal layered solid torus of type II which contains an edge of degree three not in $\cup Y_i.$ Let $T$ be such a maximal layered solid torus, and let $e$ be its unique $0$--even boundary edge. Since $\tri$ is (II,4)--free, we have $d(e)\ge 5.$ If $e$ is the longitude of a maximal layered solid torus of type IV, then $\e (T)$ is contained in one of the sets $Y_i.$ Hence assume this is not the case. We make the following observations (with the amount of detail useful for future applications):
\begin{enumerate}
\item If $e$ is not incident with another maximal layered solid torus of type II containing a degree three edge, and $T$ is almost supportive, then $\{ e\}$ is associated with $\{\e(T)\}$ and gives a cancellation.
\item If $e$ is not incident with another maximal layered solid torus of type II containing a degree three edge, and $T$ is not almost supportive, then $T$ also has an $H$--even interior edge, $e_T,$ of degree at least five. Hence associate $\{e, e_T\}$ with $\{\e(T)\},$ giving a deficit.
\item Suppose $T_1$ is the unique maximal layered solid torus of type II containing a degree three edge which meets $T$ in $e.$ 
\item[(3a)] Assume that no other maximal layered solid torus of type II contains $e.$ Then Lemma \ref{lem:two meet in edge} implies that $T, T_1$ are the only maximal layered solid tori of type II meeting in an $H$--even edge. 

If $d(e)\ge 7,$ then $\{e \}$ is associated with $\{\e(T), \e(T_1)\}$ and we have a deficit.

If $d(e)= 6$ and one of $T, T_1$ contains an $H$--even interior edge, $e'$ of degree at least five, then $\{e, e'\}$ is associated with $\{\e(T), \e(T_1)\}$ and gives a deficit. If $d(e)= 6$ and all $H$--even interior edges of $T, T_1$ are of degree four, then $\{ e\}$ is associated with $\{\e(T), \e(T_1)\}$ and gives a cancellation.

If $d(e)=5$ and both $T, T_1$ are almost supportive, then $\{e\}$ is associated with $\{\e(T), \e(T_1)\}$ and gives a gain of $+1.$ If $d(e)=5$ and one of $T, T_1$ is almost supportive and the other contains a unique $H$--even interior edge $e'$ of degree five and all others have degree four, then $\{e, e'\}$ is associated with $\{\e(T), \e(T_1)\}$ and gives a cancellation. Otherwise a deficit is associated with at most three $H$--even edges. If $d(e)=5$ and neither of $T, T_1$ is almost supportive, then there are distinct $H$--even interior edges, $e', e''$ of degrees at least five and $\{e, e', e''\}$ is associated with $\{\e(T), \e(T_1)\}$ and gives a deficit. 

\item[(3b)] Suppose the maximal layered solid torus of type II $T_2$ meets $T$ and $T_1$ in $e.$ Lemma \ref{lem:three meet in edge} implies that $T, T_1, T_2$ are the only maximal layered solid tori of type II in the triangulation. Hence $Y = \cup Y_i \bigcup \{\e(T), \e(T_1)\}.$ Since no two of the three can meet in a face, we have $d(e) \ge 6.$ It follows that there is a cancellation if all $H$--even interior edges have degree four, and a deficit associated to $e$ together with a collection of interior edges otherwise.

\item Suppose $T_1$ and $T_2$ are distinct maximal layered solid tori of type II containing degree three edges meeting $T$ in $e.$ As above, $Y = \cup Y_i \bigcup \{\e(T), \e(T_1), \e(T_2)\}$ and $d(e) \ge 6.$ Hence $e$ together with a collection of interior edges can be associated with $\{\e(T), \e(T_1), \e(T_2)\}$ giving a deficit unless one of the following cases occurs. If $d(e) = 7$ and every $H$--even interior edge of $T, T_1, T_2$ has degree four, then there is a cancellation. If $d(e)=6$ and every $H$--even interior edge of $T, T_1, T_2$ has degree four, then there is a gain of $+1.$ If $d(e)=6$ and there is precisely one $H$--even interior edge of degree five and all others have degree four, then there is a cancellation.
\end{enumerate}
It follows that the maximal gain which can be achieved is $+1;$ whence (\ref{eqn:degree three proof}) cannot be satisfied.
\end{proof}


\subsection{Proofs of the main results}

\begin{proposition}\label{cor:norm bound 2}
Let $M$ be a closed, orientable, irreducible 3--manifold with minimal triangulation $\tri.$ Suppose that $H$ is a rank--2 subgroup of $H^1(M; \Z_2).$ Then
$$|\tri| \ge 2 + \sum_{0\neq \varphi\in H} || \ \varphi \ ||.$$
Suppose further that all edge loops are coloured by $H$ and that $\tri$ is (II,4)--free. Then
$$|\tri| - \sum_{0\neq \varphi\in H} || \ \varphi \ || \ge |\tri| + \sum_{0\neq \varphi\in H} \chi(S_\varphi) \ge 2.$$
\end{proposition}
\begin{proof}
First replace $\tri$ by a (II,4)--free minimal triangulation, $\tri_0.$ Then $|\tri|=|\tri_0|.$
If $|\tri_0| + \sum \chi(S_\varphi) \le 1,$ then equation (\ref{inequ for min tri general}) in Lemma~\ref{lem:formula for degree 3 edges} and inequality (\ref{eqn:inequality degree three}) in Lemma~\ref{lem:bound for degree 3 edges} give a contradiction. Hence 
$$|\tri_0| - \sum_{0\neq \varphi\in H} || \ \varphi \ || \ge |\tri_0| + \sum_{0\neq \varphi\in H} \chi(S_\varphi) \ge 2,$$
where the first inequality follows from Lemma \ref{lem:chi bounds norm for canonical}.
\end{proof}

The first part of the above proposition directly implies Theorem \ref{thm:general bound} by definition of a minimal triangulation. The implication (2) $\Longrightarrow$  (1) in Theorem \ref{thm:really main} is a consequence of Lemma \ref{lem:quad surfaces} and Theorem \ref{thm:general bound}. The reverse is contained in the following statement:

\begin{proposition}\label{thm:gen really main}
Let $M$ be a closed, orientable, irreducible 3--manifold with arbitrary one-vertex triangulation $\tri.$ Suppose that $H$ is a rank--2 subgroup of $H^1(M; \Z_2)$ and that 
$$ |\tri|  = 2+  \sum_{0\neq \varphi\in H} || \ \varphi \ ||.$$
Then $\tri$ is minimal. Moreover, every minimal triangulation of $M$ has the property that each canonical surface dual to a non-zero element of $H$ is a $\Z_2$--taut quadrilateral surface.
\end{proposition}

\begin{proof}
Every minimal triangulation of $M$ has a single vertex. Let $\tri_0$ be a minimal triangulation which is (II,4)--free with respect to the colouring of edges by $H,$ and denote $S_i$ the dual surfaces. Then
$$2 \le |\tri_0| + \sum \chi(S_i) \le |\tri| - \sum || \ \varphi_i \ || = 2$$
by Proposition \ref{cor:norm bound 2} and Lemma \ref{lem:chi bounds norm for canonical}. This forces equality and since $ |\tri_0| \le  |\tri|$ and $\chi(S_i) \le -|| \ \varphi_i \ ||,$ we have equality for all terms. In particular, $\tri$ is minimal and each canonical surface is taut. Lemma \ref{lem:formula for degree 3 edges} implies that $A(\tri_0)+C(\tri_0) \le 1,$ and using Lemma \ref{lem:char of small cases}, this implies $E(\tri_0)=|\tri_0|.$ This gives the desired conclusion for a (II,4)--free minimal triangulation. Given an arbitrary minimal triangulation which is not (II,4)--free, one may perform a finite sequence of edge flips to arrive at a (II,4)--free minimal triangulation. The latter must have an edge of degree four which is contained in four distinct tetrahedra. Since each tetrahedron is of type V, this implies that there is a compression disc for one of the canonical surfaces, contradicting the fact that it is taut.
\end{proof}

\begin{remark}
Under the hypothesis of Proposition \ref{thm:gen really main}, if there is a rank-2 subgroup $G\neq H,$ then
$$\sum_{0\neq \varphi\in G} || \ \varphi \ || < \sum_{0\neq \varphi\in H} || \ \varphi \ ||.$$
\end{remark}


\section{The twisted layered loop triangulation}
\label{sec:triangulations}

This section shows that the twisted layered loop triangulation of $M_k=S^3/Q_{4k},$ $k\ge 1,$ satisfies the hypothesis of Theorem \ref{thm:really main} when $k$ is even. If $k$ is odd, then $H_1(M_k; \Z) = \Z_4$ and the results of this paper do not apply.


\subsection{The twisted layered loop triangulation}
\label{sec:Layered loop triangulation}

\begin{figure}[t]
\psfrag{i}{{\small $e_2$}}
\psfrag{h}{{\small $e_1$}}
\psfrag{b}{{\small $b$}}
\psfrag{t}{{\small $t$}}
\psfrag{e}{{\small $e_h$}}
\psfrag{f}{{\small $e_{h+1}$}}
\psfrag{g}{{\small $e_{h+2}$}}
\begin{center}
      \includegraphics[height=3.5cm]{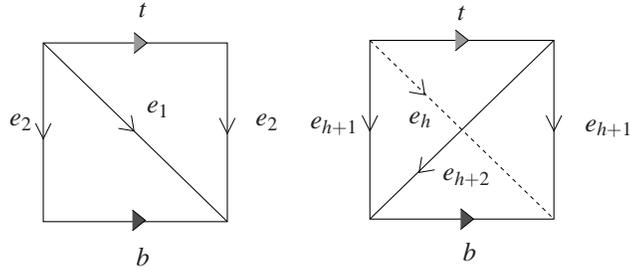}
\end{center}
    \caption{The twisted layered loop triangulation}
     \label{fig:quaternionic}
\end{figure}

Starting point is the triangulation with two faces of the annulus shown with labelling in Figure \ref{fig:quaternionic}. The edges corresponding to the two boundary components are denoted $t$ for \emph{top} and $b$ for \emph{bottom}, and oriented so that they correspond to the same element in fundamental group. The remaining two edges are $e_1$ and $e_2,$ oriented from $t$ to $b.$ Tetrahedron $\sigma_1$ is layered along $e_1,$ and the new edge denoted $e_3$ and oriented from $t$ to $b.$ The annulus is thus identified with two faces of $\sigma_1.$ Inductively, tetrahedron $\sigma_h$ is layered along edge $e_h,$ and the new edge $e_{h+2}$ is oriented from $t$ to $b.$ Assume $k$ tetrahedra have thus been attached; if $k=0$ we have an annulus, if $k=1$ a creased solid torus and if $k \ge 2$ a solid torus. Denote the resulting triangulation $C_k.$

The two free faces of tetrahedron $\sigma_k$ in $C_k$ are identified with the two free faces of tetrahedron $\sigma_1$ such that $\sigma_1$ is layered along $e_{k+1}$ with $e_1 \leftrightarrow -e_{k+1},$ $e_2 \leftrightarrow -e_{k+2}$ and $t\leftrightarrow-b.$ The result is $M_k$ (see \cite{bab}, Theorem 3.3.11, or \cite{JRT2}), and the triangulation, denoted $\widehat{C}_k,$ is termed its \emph{twisted layered loop triangulation}.


\subsection{One-sided incompressible surfaces}

Note that there is a normal, one-sided Klein bottle, $S_1,$ in $M_k,$ obtained by placing a quadrilateral in each tetrahedron dual to the edge $t.$ When $k$ is even, then there are two more embedded normal surfaces, $S_2$ and $S_3,$ just containing quadrilateral discs. One is dual to the set of all edges $e_i,$ where $i$ is even; the other is dual to the set of all edges $e_i$ where $i$ is odd. We have that $M_k \setminus S_1$ is a torus; and $M_k \setminus S_2$ and $M_k \setminus S_3$ are genus $\frac{k}{2}$ handlebodies.
Using singular homology with $\Z_2$--coefficients, the intersection pairing induces homomorphisms $\varphi_i \co \pi_1(M_k) \to \Z_2$ defined by:
\begin{itemize}
\item[] $\varphi_1[e_i] = 1$ for each $i$ and $\varphi_1[t]=0.$
\item[] $\varphi_2[e_i] = 0$ if $i$ is even, $\varphi_2[e_i] = 1$ if $i$ is odd and $\varphi_2[t]=1;$ 
\item[] $\varphi_3[e_i] = 0$ if $i$ is odd, $\varphi_3[e_i] = 1$ if $i$ is even and $\varphi_3[t]=1.$
\end{itemize}
We have $\varphi_1 + \varphi_2 = \varphi_3,$ and  $H^1(M; \Z_2) = \langle \varphi_1, \varphi_2 \rangle.$ Since $S_i$ is dual to $\varphi_i,$ it follows that we have found a representative for each non-trivial $Z_2$--cohomology class. To show that the triangulation satisfies the hypothesis of Theorem \ref{thm:really main}, we need to show that these surfaces are in fact representatives of minimal genus. This follows from the following application of work by Waldhausen \cite{wald}, Rubinstein \cite{Rubin1978} and Frohman \cite{Fro1986}.
\begin{lemma}\label{lem:one-sided in quaternionic}
If $k$ is even, then $S_1,$ $S_2$ and $S_3$ are, up to isotopy, the only connected, one-sided incompressible surfaces in $M_k.$
\end{lemma}
\begin{proof}
Let $N_1$ be a small regular neighbourhood of $S_1$ in $M_k.$ Then $N_1$ is a twisted $I$--bundle over the Klein bottle. Moreover, $N_1$ is homeomorphic to the $S^1$--bundle over the M\"obius band with orientable total space. This has a Seifert fibration over the disc with two cone points of order two. Denote the base orbifold $D,$ its cone points $p_1$ and $p_2,$ and the fibration $p\co N_1 \to D.$ Equivalently, this description of $N_1$ can be obtained by cutting $M_k=S^2(\ (2,1), (2,1), (k,1-k) \ )$ along the vertical torus which is the pre-image of the boundary of a small regular neighbourhood of the cone point labelled $(k,1-k).$

Waldhausen (\cite{wald}, Section 3) showed that the only connected, two-sided incompressible and non-boundary parallel surfaces in $N_1$ are a vertical annulus which is the pre-image of a properly embedded arc in $D$ separating $p_1$ and $p_2;$ and a horizontal annulus which is a branched double cover of $D.$ Using the methods of Frohman \cite{Fro1986}, one can similarly show that the only connected, one-sided incompressible surfaces in $N_1$ are a Klein bottle, $K,$ (which is the pre-image in $N_1$ of an arc in $D$ joining $p_1$ and $p_2$) as well as two M\"obius bands, $B_1$ and $B_2,$ (where $B_i$ is the pre-image in $N_1$ of an arc in $D$ joining $p_i$ to $\partial D$). Note that the M\"obius bands have parallel boundary curves on $\partial N_1.$

Similarly, it is known through the work of Waldhausen (\cite{wald}, Section 2) that a connected, two-sided incompressible surface in the solid torus $M_k \setminus N_1$ is either a meridian disc or a boundary-parallel disc or annulus. Rubinstein \cite{Rubin1978} showed that a one-sided incompressible surface in $M_k \setminus N_1$ has boundary a single curve which uniquely determines the surface up to isotopy.

Suppose $S$ is a connected, one-sided, incompressible surface in $M_k,$ and isotope $S$ such that it meets $\partial N_1$ minimally in a collection of pairwise disjoint curves. Then $S\cap \partial N_1$ is either empty or consists of finitely many parallel, essential curves on $\partial N_1.$

If $S\cap \partial N_1 = \emptyset,$ then $S=K,$ since the complement of $N_1$ is a solid torus. Hence assume $S\cap \partial N_1 \neq \emptyset.$ We have that $S \cap N_1$ (respectively $S \cap (M_k \setminus N_1)$) is incompressible in $N_1$ (respectively $M_k \setminus N_1$).

First assume that $S \cap N_1$ is two-sided. Then $S\cap N_1$ is a family of parallel annuli and $S\cap \partial N_1$ has an even number of components. It follows that $S \cap (M_k \setminus N_1)$ must also be a family of parallel annuli. Since an outermost annulus in the latter family is boundary compressible, there is an isotopy of $S$ reducing the number of curves in  $S\cap \partial N_1.$ This contradicts the assumption that $S$ meets $\partial N_1$ minimally.

Hence assume that $S \cap N_1$ is one-sided. It follows that $S$ is either $B_1,$ $B_2$ or $B_1 \cup B_2.$ If $S\cap N_1$ is $B_1$ or $B_2,$ then $S \cap (M_k \setminus N_1)$ has a single boundary component and hence is a uniquely determined incompressible one-sided surface, $F.$ If $S\cap N_1=B_1 \cup B_2,$ then $S \cap (M_k \setminus N_1)$ is a boundary parallel annulus, and there is an isotopy making $S\cap \partial N_1 = \emptyset,$ which contradicts the assumption that $S$ meets $\partial N_1$ minimally.

It follows that, up to isotopy, a connected, one-sided, incompressible surfaces in $M_k$ is one of $K,$ $B_1 \cup F$ and $B_2 \cup F.$ Since $H^1(M; \Z_2)$ has rank two, each of these three surfaces must be incompressible.

It follows from the construction that $K = S_1.$ We claim that, up to re-labeling, $S_2 = B_1 \cup F$ and $S_3 = B_2 \cup F.$ Notice that the boundary curve of $F$ is a regular fibre, and hence a $(k,1)$--curve on  $M_k \setminus N_1$ with respect to the standard longitude and meridian. It now follows from \cite{Rubin1978} (see also \cite{Fro1986} and \cite{JR:LT}) that 
$$\chi(S_2)=\chi(S_3)=\chi(B_1 \cup F) = \chi(B_2 \cup F).$$
This completes the proof.
\end{proof}




\address{Department of Mathematics, Oklahoma State University, Stillwater, OK 74078-1058, USA}
\email{jaco@math.okstate.edu}

\address{Department of Mathematics and Statistics, The University of Melbourne, VIC 3010, Australia} 
\email{rubin@ms.unimelb.edu.au} 

\address{Department of Mathematics and Statistics, The University of Melbourne, VIC 3010, Australia} 
\email{tillmann@ms.unimelb.edu.au} 
\Addresses
                                                      

\begin{thebibliography}{99}

\bibitem{bab} Benjamin A.\thinspace Burton: \emph{Minimal triangulations and normal surfaces}, PhD thesis, The University of Melbourne, 2003.

\bibitem{Fro1986} Charles Frohman: \emph{One-sided incompressible surfaces in Seifert fibered spaces}, Topology Appl.  23  (1986),  no. 2, 103--116.

\bibitem{JR} William Jaco and J.\thinspace Hyam Rubinstein: \emph{0--efficient triangulations of 3--manifolds}, Journal of Differential Geometry {\bf 65} (2003), no. 1, 61--168.

\bibitem{JR:LT} William Jaco and J.\thinspace Hyam Rubinstein: \emph{Layered-triangulations of 3--manifolds}, arXiv:math/0603601.

\bibitem{JRT} William Jaco, J.\thinspace Hyam Rubinstein and Stephan Tillmann: \emph{Minimal triangulations for an infinite family of lens spaces}, Journal of Topology 2 (2009) 157-180.

\bibitem{JRT2} William Jaco, J.\thinspace Hyam Rubinstein and Stephan Tillmann: \emph{Coverings and minimal triangulations of 3--manifolds}, arXiv:0903.0112v1.

\bibitem{kr} Ensil Kang and J. Hyam Rubinstein: \emph{Ideal triangulations of 3--manifolds I; spun normal surface theory}, Geometry and Topology Monographs, Vol. 7 , Proceedings of the Casson Fest, 235--265 (2004).

\bibitem{MP} Bruno Martelli and Carlo Petronio: \emph{Three-Manifolds Having Complexity At Most 9}, Experimental Mathematics 10 (2001), 207-237.
 
\bibitem{Mat1990} Sergei V.\thinspace Matveev: \emph{Complexity theory of three-dimensional manifolds}, Acta Appl. Math.  19 (1990),  no. 2, 101--130.

\bibitem{Mat2003} Sergei V.\thinspace Matveev: \emph{Complexity of three-dimensional manifolds: problems and results,} Siberian Adv. Math.  13  (2003),  no. 3, 95--103.

\bibitem{Mat2005} Sergei V.\thinspace Matveev: \emph{Recognition and tabulation of three-dimensional manifolds}, (Russian)  Dokl. Akad. Nauk  400  (2005),  no. 1, 26--28.

\bibitem{MatPer2001} Sergei V.\thinspace Matveev and Ekaterina L.\thinspace Pervova: \emph{Lower bounds for the complexity of three-dimensional manifolds}, (Russian)  Dokl. Akad. Nauk  378  (2001),  no. 2, 151--152.

\bibitem{MPV} Sergei V.\thinspace Matveev, Carlo Petronio and Andrei Vesnin: \emph{Two-sided aymptotic bounds for the complexity of some closed hyperbolic three-manifolds}, to appear in the Journal of the Australian Mathematical Society, arXiv:math/0602372v1.

\bibitem{Rubin1978} J.\thinspace Hyam Rubinstein: \emph{One-sided Heegaard splittings of 3--manifolds}, Pacific J.\thinspace Math.\thinspace 76 (1978), no. 1, 185--200.

\bibitem{Th86} W.\thinspace Thurston: \emph{A norm for the homology of $3$-manifolds.} Mem. Amer. Math. Soc. 59 (1986), no. 339, i--vi and 99--130.

\bibitem{ti} Stephan Tillmann: \emph{Normal surfaces in topologically finite 3--manifolds}, L'Enseignement Math\'ematique (2) 54 (2008), 329--380.

\bibitem{to} Jeffrey L.\thinspace Tollefson: \emph{Normal surface $Q$--theory}, Pacific J.  of Math.183 (1998), 359-374.

\bibitem{wald} Friedhelm Waldhausen: \emph{Eine Klasse von $3$-dimensionalen Mannigfaltigkeiten I, II,} Invent. Math. 3 (1967), 308--333; ibid. 4  1967 87--117.

\end{thebibliography}
\end{document}